\documentclass[12pt,a4paper,twoside]{article}
\usepackage[english]{babel}
\usepackage{amssymb}
\usepackage{inputenc}
\newcommand{\ba}{\begin{array}}
\newcommand{\ea}{\end{array}}
\usepackage{amssymb}
\begin{document}

\author{S. Albeverio $^{1},$ Sh. A. Ayupov $^{2},$ \ \ K. K.
Kudaybergenov  $^3$}
\title{\bf Description of Derivations  on Measurable Operator Algebras of   Type I}

\maketitle

\begin{abstract}

Given a type I von Neumann algebra  $M$ with a faithful normal
semi-finite trace $\tau,$ let $L(M, \tau)$ be the algebra of all
$\tau$-measurable operators affiliated with $M.$ We give a complete
description of all derivations on the algebra $L(M, \tau).$ In
particular, we prove that if $M$ is of type I$_\infty$ then every
derivation on $L(M, \tau)$ is inner.

\end{abstract}

\medskip
$^1$ Institut f\"{u}r Angewandte Mathematik, Universit\"{a}t Bonn,
Wegelerstr. 6, D-53115 Bonn (Germany); SFB 611, BiBoS; CERFIM
(Locarno); Acc. Arch. (USI), \emph{albeverio@uni-bonn.de}

$^2$ Institute of Mathematics and information  technologies,
Uzbekistan Academy of Science, F. Hodjaev str. 29, 100125, Tashkent
(Uzbekistan), e-mail: \emph{sh\_ayupov@mail.ru,
e\_ayupov@hotmail.com}

 $^{3}$ Institute of
Mathematics and information  technologies, Uzbekistan Academy of
Science, F. Hodjaev str. 29, 100125, Tashkent (Uzbekistan), e-mail:
\emph{karim2006@mail.ru}

\medskip \textbf{AMS Subject Classifications (2000): 46L57, 46L50, 46L55,
46L60}

\textbf{Key words:}  von Neumann algebras,  non commutative
integration,  measurable operator, Kaplansky-Hilbert module, type I
algebra, derivation, inner derivation.

\newpage
 \large

\section*{\center 1. Introduction}

The  present paper is devoted to  a complete description of
derivations on the   algebra of $\tau$-measurable operators $L(M,
\tau)$ affiliated with a type I von Neumann algebra $M$ and a normal
faithful semi-finite trace $\tau.$

Given a (complex) algebra $A,$ a linear operator $D:A\rightarrow A$
   is called a  \emph{derivation} if $D(xy)=D(x)y+xD(y)$ for all $x, y\in A.$
   Each element $a\in A$ generates a derivation $D_a:A\rightarrow A$
   defined as $D_a(x)=ax-xa,\,  x\in A.$ Such derivations are called  \emph{inner} derivations.

It is well known that all derivation on a von Neumann algebra are
inner and therefore are norm continuous. But the properties of
derivations on the unbounded operator algebra $L(M, \tau)$ seem to
be very far from being similar. Indeed, the results of
  \cite{Ber} and  \cite{Kus} show that in the commutative case where
 $M=L^{\infty}(\Omega, \Sigma, \mu),$ with $(\Omega, \Sigma, \mu)$
  any  non atomic measure space with a finite measure $\mu,$ the
 algebra $L(M, \tau)\cong L^{0}(\Omega, \Sigma, \mu)$ of all  complex measurable functions on
$(\Omega, \Sigma, \mu)$ admits non zero derivations. It is clear
that these derivations are  discontinuous in the  measure topology
(i. e. the  topology of convergence in measure),  and thus are non
inner. It seems that the existence of such pathological examples
deeply depends on the commutativity of the underlying algebra $M.$
Indeed, the main result of our previous paper \cite{AlAyup} states
that if $M$ is a type I von Neumann algebra, then any derivation $D$
on  $L(M, \tau),$ which is identically zero on the center $Z$ of the
von Neumann algebra $M$ (i.e. $D$ is $Z$-linear), is inner, i.e.
$D(x)=ax-xa$ for an appropriate element $a\in L(M, \tau).$

In the mentioned paper \cite{AlAyup} we have also constructed an
example of a non inner derivation on the algebra $L(M, \tau),$ where
$M$ is a homogeneous type I$_n$ algebra
$L^{\infty}(\Omega,\mu)\bar{\otimes}M_n(\mathbb{C}).$ In this case
$L(M, \tau)$ coincides with the algebra $M_n(L^{0})$ of all $n\times
n$ matrices over the algebra $L^{0}=L^{0}(\Omega, \Sigma, \mu).$
Namely, given any non zero derivation $\delta: L^{0}\rightarrow
L^{0}$ and a matrix $(\lambda_{ij})_{i, j=1}^{n}\in M_n(L^{0}),$
$\lambda_{ij}\in L^{0}, i, j=\overline{1, n},$ put
$$D_{\delta}((\lambda_{ij})_{i, j=1}^{n})=(\delta(\lambda_{ij})_{i,
j=1}^{n}).\eqno(1)$$ Then it is clear that $D_\delta$ defines a
derivation on $M_n(L^{0}),$ which coincides with $\delta$ on the
center of $M_n(L^{0}).$

In the present paper we prove that for type I von Neumann algebras
the above construction (1) gives the general form of the
pathological derivations and these only exist in type $I_{fin}$
cases, while for type $I_{\infty}$ von Neumann algebras $M$ all
derivation on $L(M, \tau)$ are inner. Moreover we prove that an
arbitrary derivation $D$ on $L(M, \tau)$ for a type I von Neumann
algebra $M,$ can be uniquely decomposed into the sum
$D_a+D_{\delta}$ where the derivation $D_a$ is inner and the
derivation $D_\delta$ is of the form (1). In such a decomposition
$\delta$ is defined uniquely, and the element $a$ is unique up to a
central summand.

\begin{center} {\bf 2. Preliminaries}
\end{center}

 Let   $(\Omega,\Sigma,\mu)$  be a measurable space and suppose
 that the measure $\mu$ has the  direct sum property, i. e. there is a family
 $\{\Omega_{i}\}_{i\in
J}\subset\Sigma,$ $0<\mu(\Omega_{i})<\infty,\,i\in J,$ such that for
any $A\in\Sigma,\,\mu(A)<\infty,$ there exist a countable subset
$J_{0 }\subset J$ and a set  $B$ with zero measure such that
$A=\bigcup\limits_{i\in J_{0}}(A\cap \Omega_{i})\cup B.$

 We
denote by  $L^{0}=L^{0}(\Omega, \Sigma, \mu)$ the algebra of all
(equivalence classes of) complex measurable functions on $(\Omega,
\Sigma, \mu)$ equipped with the topology of convergence in measure.
Then $L^{0}$ is a complete metrizable commutative regular algebra
with the unit $\textbf{1}$ given by $\textbf{1}(\omega)=1,
\,\omega\in\Omega.$

Recall that a net $\{\lambda_\alpha\}$ in $L^{0}$ $(o)$-converges to
$\lambda\in L^{0}$ if there exists a net $\{\xi_\alpha\}$ monotone
decreasing to zero such that
$|\lambda_\alpha-\lambda|\leq\xi_\alpha$ for all $\alpha.$

Denote by $\nabla$ the complete Boolean algebra of all idempotents
from $L^{0},$ i. e. $\nabla=\{\tilde{\chi}_{A}: A\in\Sigma\},$ where
$\tilde{\chi}_{A}$ is  the element from $L^{0}$ which contains the
characteristic function of the set $A.$ A \emph{partition of the
unit} in $\nabla$ is a family $(\pi_\alpha)$ of orthogonal
idempotents from $\nabla$ such that
$\sum\limits_{\alpha}\pi_\alpha=\textbf{1}.$

A complex linear space  $E$ is said to be normed  by $L^{0}$ if
there is
 a map  $\|\cdot\|:E\longrightarrow L^{0}$ such that for any  $x,y\in E, \lambda\in
\mathbb{C},$ the following conditions are fulfilled:

 1) $\|x\|\geq 0; \|x\|=0\Longleftrightarrow x=0$;

 2) $\|\lambda x\|=|\lambda|\|x\|$;

 3)  $\|x+y\|\leq\|x\|+\|y\|$.

 The  pair  $(E,\|\cdot\|)$ is called a lattice-normed  space over $L^{0}.$
A lattice-normed  space  $E$ is called  $d$-decomposable, if for any
$x\in E$ with  $\|x\|=\lambda_{1}+\lambda_{2},$  $ \lambda_{1},
\lambda_{2} \in L^{0},\,\lambda_{1}\lambda_{2}=0,$ there exist
$x_{1}, x_{2}\in E$ such that $x=x_{1}+x_{2}$ and
$\|x_{i}\|=\lambda_{i},\,i=1,2$. A net $(x_{\alpha})$ in $E$ is said
to be $(bo)$-convergent to $x\in E$, if the net
$\{\|x_{\alpha}-x\|\}$ $(o)$-converges to zero in
 $L^{0}.$
 A lattice-normed  space $E$ which is $d$-decomposable and  complete with respect to the  $(bo)$-convergence
  is called a
 \emph{Banach-Kantorovich space}.

 It is known  that every Banach-Kantorovich space $E$ over
 $L^{0}$ is a module over $L^{0}$ and  $\|\lambda x\|=|\lambda|\|x\|$
 for all $\lambda\in L^{0},\, x\in E$ (see \cite{Kusr}).

 Let     $E$ be a Banach-Kantorovich space  over $L^{0}.$
              If $(u_{\alpha})_{\alpha\in A}$ is a net in $E$ and
$(\pi_{\alpha})_{\alpha\in A}$
              is a partition of the unit in $\nabla$, then the
              series
      $\sum\limits_{\alpha}\pi_{\alpha}u_{\alpha}$ $(bo)$-converges in
      $E$ and  its sum is called the \emph{mixing} of
       $(u_{\alpha})_{\alpha\in A}$ with respect to  $(\pi_{\alpha})_{\alpha\in A}.$ We denote this sum
             by $\textmd {mix}(\pi_{\alpha}u_{\alpha})$.
      A subset $K\subset E$ is called \emph{cyclic}, if
       $\textmd {mix}(\pi_{\alpha}u_{\alpha})\in K$ for each net $(u_{\alpha})_{\alpha\in A}\subset K$ and
       for any partition  of the unit $(\pi_{\alpha})_{\alpha\in A}$ in $\nabla.$

Let $K$ be a cyclic set. We say that a  map $T:K\rightarrow K$
commutes with the mixing operation, if
$$T(\sum\limits_{\alpha}\pi_{\alpha}x_{\alpha})=
\pi_{\alpha}T(\sum\limits_{\alpha}x_{\alpha})$$ for all
$(x_\alpha)\subset L^{0}$ and   any partition $\{\pi_{\alpha}\}$ of
the unit in $\nabla.$

Let  $\mathcal{K}$ be a module over $L^{0}$. A map $\langle
\cdot,\cdot\rangle:\mathcal{K}\times \mathcal{K}\rightarrow L^{0}$
is called  an $L^{0}$-valued inner product, if for all $x,y,z\in
\mathcal{K},\,\lambda\in L^{0},$ the following conditions are
fulfilled:

1) $\langle x,x\rangle\geq0$; $\langle x,x\rangle=0\Leftrightarrow
x=0$;

2) $\langle x,y\rangle=\overline{\langle y,x\rangle}$;

3) $\langle \lambda x,y\rangle=\lambda\langle x,y\rangle$;

4)  $\langle x+y,z\rangle=\langle x,z\rangle+\langle y,z\rangle$.

If $\langle \cdot,\cdot\rangle:\mathcal{K}\times
\mathcal{K}\rightarrow L^{0}$ is an $L^{0}$-valued inner product,
then $ \|x\|=\sqrt{\langle x,x \rangle} $ defines  an $L^{0}$-valued
norm on $\mathcal{K}.$ The  pair $(\mathcal{K},\langle
\cdot,\cdot\rangle)$ is called a \emph{Kaplansky-Hilbert module}
over $L^{0},$ if $(\mathcal{K},\|\cdot\|)$ is a Banach-Kantorovich
space over  $L^{0}$ (see \cite{Kusr}).

Let  $X$ be a Banach space.  A map  $s:\Omega\rightarrow X$ is said
to be simple, if
$s(\omega)=\sum\limits_{k=1}^{n}\chi_{A_{k}}(\omega)c_k,$ where
$A_k\in\Sigma, A_i\cap A_j=\emptyset, \,i\neq j,\, \,c_k\in
X,\,k=\overline{1, n},\, n\in\mathbb{N}.$ A map $u:\Omega\rightarrow
X$ is said to be measurable, if there is a sequence  $(s_n)$ of
simple maps such that $\|s_n(\omega)-u(\omega)\|\rightarrow0$
almost everywhere on any $A\in\sum$ with $\mu(A)<\infty.$

Let $\mathcal{L}(\Omega, X)$ be the set of all measurable maps from
$\Omega$ into $X,$ and let $L^{0}(\Omega, X)$ denote the space of
all equivalence classes  with respect to the equality almost
everywhere. Denote by $\hat{u}$ the equivalence class from
$L^{0}(\Omega, X)$ which contains the measurable map $u\in
\mathcal{L}(\Omega, X).$ Further we shall identify the element $u\in
\mathcal{L}(\Omega, X)$ and the class $\hat{u}.$ Note that the
function  $\omega \rightarrow \|u(\omega)\|$
       is measurable for any $u\in \mathcal{L}(\Omega, X).$ The equivalence class containing the function
              $\|u(\omega)\|$ is denoted by
       $\|\hat{u}\|$. For  $\hat{u}, \hat{v}\in L^{0}(\Omega, X), \lambda\in L^{0}$ put
$\hat{u}+\hat{v}=\widehat{u(\omega)+v(\omega)},
\lambda\hat{u}=\widehat{\lambda(\omega) u(\omega)}.$

It is known  \cite{Kusr} that $(L^{0}(\Omega, X), \|\cdot\|)$ is a
Banach-Kantorovich space over $L^{0}.$

Put $L^{\infty}(\Omega, X)=\{x\in L^{0}(\Omega, X):\|x\|\in
L^{\infty}(\Omega)\}.$  Then $L^{\infty}(\Omega, X)$ is a Banach
space with respect to  the norm
$\|x\|_{\infty}=\|\|x\|\|_{L^{\infty}(\Omega)}.$

If $H$ is a Hilbert space, then $L^{0}(\Omega, H)$ can be equipped
with an  $L^{0}$-valued inner product $\langle x,
y\rangle=\widehat{(x(\omega), y(\omega))},$ where $(\cdot, \cdot)$
is the inner product on $H.$

Then $(L^{0}(\Omega, H),\langle \cdot,\cdot\rangle)$ is a
Kaplansky-Hilbert module over $L^{0}.$

Let  $E$ be a  Banach-Kantorovich space over $L^{0}.$
 An operator   $T: E\rightarrow E$
 is   $L^{0}$-linear if  $T(\lambda_1 x_1 +\lambda_2 x_2)=\lambda_1 T(x_1)+\lambda_2
 T(x_2)$
for all  $\lambda_1, \lambda_2\in L^{0}, \, x_1, x_2\in E.$ An
$L^{0}$-linear operator
 $T:E\rightarrow E$ is  $L^{0}$-bounded if there exists an element $c\in L^{0}$
 such that $\|T(x)\|\leq
 c\|x\|$ for any $x\in E.$  For an  $L^{0}$-bounded $L^{0}$-linear
operator $T$  we put $\|T\|=\sup\{\|T(x)\|:\|x\|\leq \textbf{1} \}.$

Let  $B(H)$ be the algebra of all bounded linear operators on a
Hilbert space  $H$ and let  $M$ be a von Neumann algebra in $B(H)$
with a faithful normal semi-finite trace $\tau.$ Denote by
$\mathcal{P}(M)$ the lattice of projections in  $M.$

A linear subspace  $\mathcal{D}$ in  $H$ is said to be affiliated
with  $M$ (denoted as  $\mathcal{D}\eta M$), if
$u(\mathcal{D})\subset \mathcal{D}$ for any unitary operator $u$
from the commutant $$M'=\{y\in B(H):xy=yx, \,\forall x\in M\}$$ of
the algebra $M.$

A linear operator  $x$ on  $H$ with the domain  $\mathcal{D}(x)$ is
said to be affiliated with  $M$ (denoted as  $x\eta M$) if
$u(\mathcal{D}(x))\subset \mathcal{D}(x)$ and $ux(\xi)=xu(\xi)$ for
all $u\in M',$ $\xi\in \mathcal{D}(x).$

A linear subspace  $\mathcal{D}$ in $H$ is called  $\tau$-dense, if

1) $\mathcal{D}\eta M;$

2) given any  $\varepsilon>0$ there exists a projection
$p\in\mathcal{P}(M)$ such that  $p(H)\subset \mathcal{D}$ and
$\tau(p^{\perp})\leq\varepsilon.$

A closed linear operator  $x$ is said to be $\tau$-\emph{measurable}
(or totally measurable) with respect to the von Neumann algebra $M,$
if $x\eta M$ and $\mathcal{D}(x)$ is  $\tau$-dense in $H.$

 Denote by $L(M,
\tau)$ the set of all $\tau$-measurable operators with respect to
 $M.$ Let $\|\cdot\|_{M}$ stand for the uniform norm in $M.$ The
\emph{measure topology,} $t_{\tau}$ in $L(M, \tau)$ is the one given
by the following system of neighborhoods of zero:
$$V(\varepsilon, \delta)=\{x\in L(M, \tau): \exists e\in\mathcal{P}(M), \tau(e^{\perp})\leq\delta, xe\in
M,  \|xe\|_{M}\leq\varepsilon\},$$ where $\varepsilon>0, \delta>0.$

It is known \cite{Nel} that $L(M, \tau)$ equipped with the measure
topology is a complete metrizable topological $\ast$-algebra.

 Let
$L^{\infty}(\Omega, \mu)\bar{\otimes}B(H)$ be the tensor product of
von  the Neumann algebras $L^{\infty}(\Omega, \mu)$ and $B(H),$ with
the trace $\tau=\mu\otimes Tr,$ where $Tr$ is the canonical trace
for operators in $B(H)$  (with its natural domain).

Denote by $L^{0}_{t}(\Omega, B(H))$ the space of equivalence classes
of point-wise Bochner measurable operator-valued  maps from $\Omega$
into $B(H)$ (see \cite{Tak}).

 Given $\hat{u}, \hat{v}\in L^{0}_{t}(\Omega,
B(H))$ put $\hat{u}\hat{v}=\widehat{u(\omega) v(\omega)},
\hat{u}^{\ast}=\widehat{u(\omega)^{\ast}}.$

Define $$L^{\infty}_{t}(\Omega, B(H))=\{x\in L^{0}_{t}(\Omega,
B(H)):\|x\|\in L^{\infty}(\Omega)\}.$$ The space
$(L^{\infty}_{t}(\Omega, B(H)), \|\cdot\|_{\infty})$ is a Banach
$\ast$-algebra.

 It is known  \cite{Tak} that the algebra  $L^{\infty}(\Omega, \mu)\bar{\otimes}B(H)$
is $\ast$-isomorphic with the algebra  $L^{\infty}_{t}(\Omega,
B(H)).$

Note also that
$$\tau(x)=\int\limits_{\Omega}Tr(x(\omega))\,d\mu(\omega).$$

Further we shall identify the algebra  $L^{\infty}(\Omega,
\mu)\bar{\otimes}B(H)$  with the algebra $L^{\infty}_{t}(\Omega,
B(H)).$

Denote by $B(L^{0}(\Omega, H))$ the algebra of all  $L^{0}$-linear
and $L^{0}$-bounded   operators on $L^{0}(\Omega, H).$

Given any $f\in L^{0}_{t}(\Omega, B(H))$  consider the element
$\Psi(f)$ from $B(L^{0}(\Omega, H))$ defined by
 $$\Psi(f)(x)=\widehat{f(\omega)(x(\omega))},\quad x\in L^{0}(\Omega, H).$$

Then the correspondence  $f\rightarrow \Psi(f)$ gives an isometric
$\ast$-isomorphism between the algebras  $L^{0}_{t}(\Omega, B(H))$
and $B(L^{0}(\Omega, H))$  (see \cite{Kusr}).

It is known \cite{AlAyup}, that the algebra $L(L^{\infty}(\Omega,
\mu)\bar{\otimes}B(H), \tau)$ of all $\tau$-measurable operators
with respect to $L^{\infty}(\Omega, \mu)\bar{\otimes}B(H)$ is
$\ast$-isomorphic with the algebra $L^{0}_{t}(\Omega, B(H)).$

Therefore one has the following relations for the algebras mentioned
above:
$$
L(L^{\infty}(\Omega)\bar{\otimes} B(H)),\tau) \cong
L^{0}_{t}(\Omega, B(H))\cong  B(L^{0}(\Omega,H))\eqno(2)
$$
($\cong$ standing for $\ast$-isomorphic).

The above isomorphisms enable us to obtain the following necessary
and sufficient condition for a derivation on the algebra $L(M,
\tau)$ to be inner.

\textbf{Theorem 2.1} \cite{AlAyup}. \emph{Let $M$ be a type I von
Neumann algebra with the center  $Z.$ A derivation  $D$ on the
algebra $L(M, \tau)$ is inner if and only if it is $Z$-linear, or
equivalently, it is identically zero on $Z.$}

\begin{center} {\bf 3. Main results}
\end{center}

Let  $A$ be an algebra with the center  $Z$ and let  $D:A\rightarrow
A$ be a derivation. Given any  $x\in A$ and a central element $z\in
Z$  we have
$$D(zx)=D(z)x+zD(x)$$
and
$$D(xz)=D(x)z+xD(z).$$
Since  $zx=xz$ and $zD(x)=D(x)z,$ it follows that  $D(z)x=xD(z)$ for
any $x\in A.$ This means that  $D(z)\in Z,$ i.e. $D(Z)\subseteq Z.$
Therefore  given any derivation $D$ on the algebra $A$ we can
consider it restriction $\delta$ onto the center $Z:$
$$\delta:z\rightarrow D(z),\,z\in Z.$$

This simple but important remark is crucial in our further
considerations.

Let $M$ be a homogeneous von Neumann algebra of type I. Then
\cite{Tak} $M$ is isomorphic to the tensor product
$L^{\infty}(\Omega)\bar{\otimes} B(H).$

First, let us consider the case  $\dim H=n<\infty.$ In this case
$B(H)$ coincides with the algebra of $M_n(\mathbb{C})$ of $n\times
n$ complex matrices, and from (2) we have that  $L(M, \tau)\cong
B(L^{0}(\Omega, H))$ is isomorphic with the algebra  $M_n(L^{0})$ of
all  $n\times n$ matrices over the algebra $L^{0}.$

In the papers \cite{Ber},  \cite{Kus} the existence of non zero
derivations on $L^{0}$ has been prove the  in the case of a non
atomic measure space $(\Omega,\Sigma,\mu).$ Given any derivation
$\delta:L^{0}\rightarrow L^{0}$ consider the elementwise  derivation
$D_\delta$ on $M_n(L^{0})$ defined as (1):
$$D_{\delta}((\lambda_{ij})_{i, j=1}^{n})=(\delta(\lambda_{ij})_{i, j=1}^{n}),
$$
where $(\lambda_{ij})_{i, j=1}^{n}\in M_n(L^{0}).$

A straightforward calculation shows that  $D_\delta$ is indeed  a
derivation on  $M_n(L^{0})$ and its restriction onto the center of
$M_n(L^{0})$ coincides with $\delta$ (recall that the center of
$M_n(L^{0})$ is isomorphic with $L^{0}$).

\textbf{Lemma 3.1.} \emph{Any derivation  $D$ on the algebra
$M_n(L^{0})$ admits a unique decomposition
$$D=D_a+D_\delta,$$
where $D_a$ is an inner derivation and  $D_\delta$ is a derivation
of the form}  (1).

Proof. Given a derivation  $D$ on $M_n(L^{0}),$ consider its
restriction $\delta$ onto the center  $Z=L^{0}$ and extend it to the
whole $M_n(L^{0})$ by the form (1) as  $D_\delta.$ Put
$D_1=D-D_\delta.$ Then given any  $z\in Z$ we have
$$D_1(z)=D(z)-D_\delta(z)=D(z)-D(z)=0,$$
i.e. $D_1$ is identically zero on $Z$ and therefore it is
$Z$-linear. Theorem 2.1 implies that $D_1$ is inner, i.e.  $D_1=D_a$
for an appropriate $a\in M_n(L^{0}).$ Therefore $D=D_a+D_\delta.$

Now suppose that
$$D=D_{a_{1}}+D_{\delta_{1}}=D_{a_{2}}+D_{\delta_{2}}.$$ Then
$D_{a_{1}}-D_{a_{2}}=D_{\delta_{2}}-D_{\delta_{1}}.$ Since
$D_{a_{1}}-D_{a_{2}}$ is identically zero on the center of
$M_n(L^{0}),$ then $D_{\delta_{2}}-D_{\delta_{1}}$ also is
identically zero on the center of $M_n(L^{0}).$ Thus
 $\delta_{1}=\delta_{2}$ and hence $D_{a_{1}}=D_{a_{2}}.$ The
proof is complete. $\blacksquare$

In order to consider the case of homogeneous type I von Neumann
algebra $L^{\infty}(\Omega)\bar{\otimes} B(H)$ with $\dim H=\infty,$
we need some auxiliary results.

 \textbf{Lemma 3.2.} \emph{Any derivation  $\delta$ on the algebra
 $L^{0}$ commutes with the mixing operation on nets in $L^{0}.$}

Proof. Consider a net  $\{\lambda_\alpha\}$ in $L^{0}$ and a
partition of the unit  $\{\pi_{\alpha}\}$ in $\nabla\subset L^{0}.$
Since  $\delta(\pi)=0$ for any idempotent  $\pi\in\nabla,$ we have
$\delta(\pi_\alpha)=0$ for all $\alpha$ and thus
$\delta(\pi_{\alpha}\lambda)=\pi_{\alpha}\delta(\lambda)$ for any
$\lambda\in L^{0}.$ Therefore for each $\pi_{\alpha_{0}}$ from the
given partition of the unit we have
$$\pi_{\alpha_{0}}\delta(\sum\limits_{\alpha}\pi_{\alpha}\lambda_{\alpha})=
\delta(\pi_{\alpha_{0}}\sum\limits_{\alpha}\pi_{\alpha}\lambda_{\alpha})=
\delta(\pi_{\alpha_{0}}\lambda_{\alpha_{0}})=
\pi_{\alpha_{0}}\delta(\lambda_{\alpha_0}).$$ By taking the sum over
all $\alpha_0$
 we obtain $$\delta(\sum\limits_{\alpha}\pi_{\alpha}\lambda_{\alpha})=
\sum\limits_{\alpha}\pi_{\alpha}\delta(\lambda_{\alpha}).$$ The
proof is complete. $\blacksquare$

\textbf{Lemma 3.3.}\emph{Given any non trivial derivation
$\delta:L^{0}\rightarrow L^{0}$ there exist a sequence
$\{\lambda_n\}_{n=1}^{\infty}$ in $L^{0}$ with} $|\lambda_n|\leq
\textbf{1},\,n\in \mathbb{N},$ \emph{and an idempotent $\pi\in
\nabla,\,\pi\neq 0$ such that
$$|\delta(\lambda_n)|\geq n\pi$$
for all $n\in \mathbb{N}.$}

Proof. Suppose that the set  $\{\delta(\lambda): \lambda\in L^{0},
|\lambda|\leq\textbf{1}\}$ is order bounded in  $L^{0}.$  Then
$\delta$ maps any uniformly convergent sequence in
$L^{\infty}(\Omega)$ to an $(o)$-convergent sequence in $L^{0}.$ The
algebra $L^{\infty}(\Omega)$ coincides with the uniform closure of
the linear span of idempotents from $\nabla.$ Since $\delta$ is
identically zero on $\nabla$ it follows that $\delta\equiv 0$ on
$L^{\infty}(\Omega).$ Since $\delta$ commutes with the mixing
operation and any element $\lambda\in L^{0}$ can be represented as
$\lambda=\sum\limits_{\alpha}\pi_{\alpha}\lambda_{\alpha},$ where
$\{\lambda_\alpha\}\subset L^{\infty}(\Omega, \mu),$ and
$\{\pi_{\alpha}\}$ is a partition of unit in  $\nabla,$ we have
$\delta(\lambda)=\delta(\sum\limits_{\alpha}\pi_{\alpha}\lambda_{\alpha})=
\sum\limits_{\alpha}\pi_{\alpha}\delta(\lambda_{\alpha})=0,$ i.e.
$\delta\equiv 0$ on $L^{0}.$ This contradiction shows that the set
$\{\delta(\lambda): \lambda\in L^{0}, |\lambda|\leq\textbf{1}\}$ is
not order bounded in  $L^{0}.$ Further, since  $\delta$ commutes
with the mixing operations and the set  $\{\lambda: \lambda\in
L^{0}, |\lambda|\leq\textbf{1}\}$ is cyclic, the set
$\{\delta(\lambda): \lambda\in L^{0}, |\lambda|\leq\textbf{1}\}$ is
also cyclic.
 By \cite[Proposition 3]{Kud} there exist a sequence
  $\{\lambda_n\}_{n=1}^{\infty}$ in $L^{0}$ and an idempotent
  $\pi\in \nabla,\,\pi\neq 0,$ such that $|\delta(\lambda_n)|\geq n\pi,\,n\in \mathbb{N}.$
  The proof is complete. $\blacksquare$

Now we are in position to consider derivations on measurable
operators for homogeneous von Neumann algebras of type I$_{\infty}.$

\textbf{Theorem 3.4.} \emph{If $\dim H=\infty,$ then any derivation
on the algebra  $L(L^{\infty}(\Omega)\bar{\otimes}B(H),\tau)$ is
inner.}

Proof. First let us suppose that  $H$ is separable and let
$\{\varphi_n\}_{n=1}^{\infty}$ be an  orthonormal basis in $H.$
Given any  $n\in \mathbb{N}$ put $e_n=\textbf{1}\otimes\varphi_n,$
i.e. $e_n$ is the equivalence class in $L^{0}(\Omega, H)$ which
contains the constant mapping $\omega\rightarrow \varphi_n.$

Then $\{e_n\}_{n=1}^{\infty}$ is a $\nabla$-orthonormal basis in
$L^{0}(\Omega, H),$ i.e. $\langle e_i,
e_j\rangle=\delta_{ij}\textbf{1},$ where  $\delta_{ij}$ is the
Kronecker  symbol, and any element $\xi\in L^{0}(\Omega, H)$ has the
form
$$\xi=\sum\limits_{k=1}^{\infty}\alpha_k e_k,$$
where $\alpha_k\in L^{0},\,
\sum\limits_{k=1}^{\infty}|\alpha_k|^{2}\in L^{0}.$

For each  $n\in \mathbb{N}$ consider the orthogonal projection
$p_n$ onto the submodule $\{\alpha e_n: \alpha\in L^{0}\}\subset
L^{0}(\Omega, H),$ i.e.
$$p_n(\sum\limits_{k=1}^{\infty}\alpha_k e_k)=\alpha_n e_n.$$

For any order bounded sequence  $\Lambda=\{\lambda_k\}$ in $L^{0}$
define an operator  $x_\Lambda$ on $B(L^{0}(\Omega, H))$ putting
$$x_\Lambda(\sum\limits_{k=1}^{\infty}\alpha_k e_k)=
\sum\limits_{k=1}^{\infty}\lambda_k \alpha_k e_k.$$ Then
$$x_\Lambda p_n=p_n x_\Lambda=\lambda_n p_n.\eqno (3)$$

Let  $D$ be a derivation on $B(L^{0}(\Omega, H)),$ and let $\delta$
be its restriction onto the center of  $B(L^{0}(\Omega, H)),$
identified with $L^{0}.$

Take any $\lambda\in L^{0}$ and $n\in \mathbb{N}.$ From the identity
$$D(\lambda p_n)=D(\lambda)p_n+\lambda D(p_n)$$
 multiplying by $p_n$ on both sides we obtain
$$p_nD(\lambda p_n)p_n=p_n D(\lambda)p_n+\lambda p_n D(p_n)p_n.$$
Since  $p_n$ is a projection, one has that  $p_n D(p_n)p_n=0,$ and
since $D(\lambda)=\delta(\lambda)\in Z,$ we have
$$p_nD(\lambda p_n)p_n=\delta(\lambda)p_n.\eqno (4)$$
Now from the identity
$$D(x_\Lambda p_n)=D(x_\Lambda)p_n+x_\Lambda D(p_n),$$
in view of (3)  one has similarly
$$p_nD(\lambda_n p_n)p_n=p_n
D(x_\Lambda)p_n+\lambda p_n D(p_n)p_n,$$ i.e.
$$p_n D(\lambda_n p_n)p_n=p_nD(x_\Lambda)p_n.\eqno (5)$$
Now (4) and (5) imply
$$p_nD(x_\Lambda)p_n=\delta(\lambda_n)p_n.$$
Further we have
$$\|p_nD(x_\Lambda)p_n\|\leq\|p_n\|\|D(x_\Lambda)\|\|p_n\|=
\|D(x_\Lambda)\|$$ and
$$\|\delta(\lambda_n)p_n\|=|\delta(\lambda_n)|.$$
Therefore
$$\|D(x_\Lambda)\|\geq|\delta(\lambda_n)|$$
for any bounded sequence $\Lambda=\{\lambda_n\}$ in $L^{0}.$

If we suppose that $\delta\neq 0$ then  by Lemma 3.3 there exist a
bounded sequence $\Lambda=\{\lambda_n\}$ in $L^{0}$ and an
idempotent $\pi\in \nabla,\,\pi\neq0,$ such that
$$|\delta(\lambda_n)|\geq n\pi $$
for any $n\in \mathbb{N}.$ Thus
$$\|D(x_\Lambda)\|\geq n\pi$$
for all $n\in \mathbb{N},$ i.e. $\pi=0$ -- which is  a
contradiction. Therefore $\delta\equiv 0,$ i.e. $D$ is identically
zero on the center of $B(L^{0}(\Omega, H)).$ By Theorem 2.1 $D$ is
inner.

Now suppose that $H$ is not necessarily separable and take a
projection $p$ in $B(H)$ such that $p(H)$ is a separable infinite
dimensional Hilbert space. Put $e=\textbf{1}\otimes p.$ Define a
derivation $D_e$ on $eB(L^{0}(\Omega, H))e$ by putting
$$D_e(x)=eD(x)e,\quad x\in eB(L^{0}(\Omega, H))e.$$
Since the algebra  $eB(L^{0}(\Omega, H))e$ is isomorphic with
$B(L^{0}(\Omega, p(H)))$ we have from the above separable case that
 $D_e$ is zero on the center of $eB(L^{0}(\Omega, H))e=L^{0}e.$  Therefore one has
 $0=eD(\lambda e)e=eD(\lambda)e+\lambda eD(e)e=eD(\lambda)e=D(\lambda)e,$ i.e.
 $D(\lambda)=0$ for all $\lambda\in L^{0}.$ Thus $D$ is identically zero on the center of
 $B(L^{0}(\Omega, H)).$ By Theorem 2.1 $D$ is inner. The proof is complete. $\blacksquare$

Now let us consider the general case of type I von Neumann algebras.

Recall that a von Neumann algebra $M$ is of \emph{type I} if it is
isomorphic to a von Neumann algebra with an abelian commutant, or,
equivalently $M$ admits a faithful abelian projection.

It is well-known \cite{Tak} that if   $M$ is a type I von Neumann
algebra
 then there is a unique (cardinal-indexed) orthogonal
family of projections $(q_{\alpha})_{\alpha\in
I}\subset\mathcal{P}(M)$ with $\sum\limits_{\alpha\in
I}q_{\alpha}=\textbf{1}$ such that  $q_{\alpha}M$ is a homogeneous
type I$_\alpha$ von Neumann algebra and it is isomorphic to the
tensor product of an abelian von Neumann algebra
$L^{\infty}(\Omega_{\alpha}, \mu_{\alpha})$ and $B(H_{\alpha})$ with
$\dim H_{\alpha}=\alpha,$ i. e.
$$M\cong\sum\limits_{\alpha}^{\oplus}L^{\infty}(\Omega_{\alpha}, \mu_{\alpha})\bar{\otimes}B(H_{\alpha}).$$

Consider the faithful normal semi-finite trace $\tau$ on $M,$
defined by
$$\tau(x)=\sum\limits_{\alpha}\tau_{\alpha}(x_{\alpha}),\quad x=(x_{\alpha})\in M,\,x\geq0,$$
where $\tau_{\alpha}=\mu_{\alpha}\otimes Tr_{\alpha}.$

Let

 $$\prod\limits_{\alpha}L(L^{\infty}(\Omega_{\alpha},
\mu_{\alpha})\bar{\otimes}B(H_{\alpha}), \tau_{\alpha})$$ be the
topological (Tychonoff) product of the spaces
$L(L^{\infty}(\Omega_{\alpha},
\mu_{\alpha})\bar{\otimes}B(H_{\alpha}), \tau_{\alpha})$ equipped
with corresponding measure topologies.

  Then (see \cite{ChilLit}) we have the topological embedding
$$L(M, \tau)\subset\prod\limits_{\alpha}L(L^{\infty}(\Omega_{\alpha},
\mu_{\alpha})\bar{\otimes}B(H_{\alpha}), \tau_{\alpha}).$$

It should be noted that the above topological imbedding is strict in
general (see \cite{ChilLit}).

Now let $\delta$  be a derivation on the center $Z$ of $L(M, \tau).$
Since $q_\alpha Z\cong L^{0}(\Omega_\alpha)$ for each $\alpha,$ and
$\delta$ maps each $q_\alpha Z$ into $q_\alpha Z,$ it follows that
$\delta$ induces a derivation $\delta_\alpha$ on each
$L^{0}(\Omega_\alpha).$

Put  $F=\{\alpha\in I: \dim H_{\alpha}=\alpha<\infty\}.$  Let
$D_{\delta_\alpha}$ $(\alpha\in F)$ be the derivation on the matrix
algebra $q_\alpha L(M, \tau)\cong M_\alpha(L^{0}(\Omega_\alpha))$
constructed  by formula (1). For $\alpha\in I\setminus F$ put
$D_{\delta_\alpha}\equiv 0.$ Now define a derivation $D_\delta$ on
$L(M, \tau)$ by putting
$$D_\delta(x)=(D_{\delta_\alpha}(x_\alpha)),\quad x=(x_\alpha)\in L(M, \tau).\eqno (6)$$

From Lemma 3.1 and Theorem 3.4 we obtain the following main result
of the paper.

\textbf{Theorem 3.5.} \emph{Let $M$ be a type I von Neumann algebra.
Any derivation  $D$ on the algebra  $L(M,\tau)$ can be decomposed in
a unique way as
$$D=D_a+D_\delta,$$
where  $D_a$ is an inner derivation and $D_\delta$ is a derivation
of the form} (6).

\textbf{Corollary 3.6.} \emph{Let $M$ be a  von Neumann algebra of
type I$_{\infty}.$  Then any  derivation  on the algebra $L(M,\tau)$
is inner.}

\vspace{1cm}

\textbf{Acknowledgments.} \emph{The second and third named authors
would like to acknowledge the hospitality of the $\,$ "Institut
f\"{u}r Angewandte Mathematik",$\,$ Universit\"{a}t Bonn (Germany).
This work is supported in part by the DFG 436 USB 113/10/0-1 project
(Germany) and the Fundamental Research  Foundation of the Uzbekistan
Academy of Sciences.}

\end{document}